%The volume of hyperbolic alternating link complements

\input psfig
\input amssym.def
\input amssym
\magnification=1100
\baselineskip = 0.25truein
\lineskiplimit = 0.01truein
\lineskip = 0.01truein
\vsize = 8.5truein
\voffset = 0.2truein
\parskip = 0.10truein
\parindent = 0.3truein
\settabs 12 \columns
\hsize = 5.4truein
\hoffset = 0.4truein

\setbox\strutbox=\hbox{%
\vrule height .708\baselineskip
depth .292\baselineskip
width 0pt}
\font\caps=cmcsc10

\def\sqr#1#2{{\vcenter{\vbox{\hrule height.#2pt
\hbox{\vrule width.#2pt height#1pt \kern#1pt
\vrule width.#2pt}
\hrule height.#2pt}}}}
\def\square{\mathchoice\sqr46\sqr46\sqr{3.1}6\sqr{2.3}4}
\def\leaderfill{\leaders\hbox to 1em{\hss.\hss}\hfill}
\font\bigtenrm=cmr10 scaled 1400
\tenrm

\centerline{\bigtenrm THE VOLUME OF HYPERBOLIC}
\centerline{\bigtenrm ALTERNATING LINK COMPLEMENTS}
\tenrm
\vskip 14pt
\centerline{MARC LACKENBY}
\centerline{\caps with an appendix by Ian Agol and Dylan Thurston}
\vskip 18pt
\centerline {\caps 1. Introduction}
\vskip 6pt

A major goal of knot theory is to relate the geometric structure
of a knot complement to the knot's more basic topological properties. 
In this paper, we will do this for hyperbolic
alternating knots and links, by showing that the link's most
fundamental geometric invariant - its volume - 
can be estimated directly from its alternating diagram.

A {\sl bigon region} in a link diagram is a complementary region
of the link projection having two crossings in its boundary.
A {\sl twist} is {\sl either} a connected collection of bigon regions
arranged in a row, which is maximal in the sense that it
is not part of a longer row of bigons, {\sl or} a single crossing
adjacent to no bigon regions. The {\sl twist
number} $t(D)$ of a diagram $D$ is its number of twists.
(See Figure 1.) Recall that a diagram is {\sl prime} if any simple
closed curve in the diagram that intersects the link
projection transversely in two points disjoint from the crossings 
bounds a disc that contains no
crossings. Menasco proved [5] that a link with a connected prime
alternating diagram, other 
than the standard diagram of the $(2,n)$-torus link, is hyperbolic.
Our main theorem is the following rather surprising result,
which asserts that the link complement's hyperbolic volume
is, up to a bounded factor, simply the diagram's twist number.

\noindent {\bf Theorem 1.} {\sl Let $D$ be a prime alternating
diagram of a hyperbolic link $K$ in $S^3$. Then
$$v_3(t(D)-2)/2 \leq {\rm Volume}(S^3 - K) < v_3 (16 t(D) - 16),$$
where $v_3 (\approx 1.01494)$ is the volume of a regular hyperbolic 
ideal 3-simplex.}

The upper bound on volume actually applies to any diagram
of a hyperbolic link, not just an alternating one. The lower
bound on volume can be improved to $v_3(t(D) - 2)$ if,
in addition, $D$ is `twist-reduced'. We will define this 
term later in the paper and show that any prime alternating link has a 
twist-reduced prime alternating diagram.

Shortly after this paper was distributed, Dylan Thurston and Ian
Agol improved the upper bound in Theorem 1 to $10 v_3(t(D)-1)$.
Moreover, they showed that this new upper bound is asymptotically
sharp, in that there is a sequence of links $K_i$ with prime
alternating diagrams $D_i$ such that ${\rm Volume}(K_i) / 10 v_3 \ t(D_i)
\rightarrow 1$ as $i \rightarrow \infty$. Their results are given
in an appendix to this paper.

\vskip 18pt
\centerline{\psfig{figure=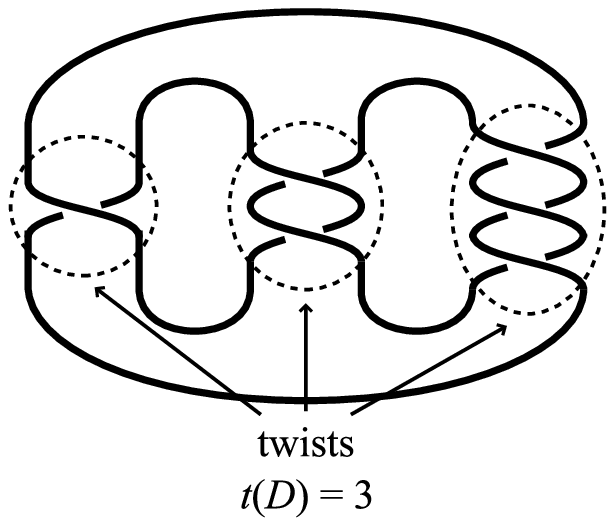,width=2in}}
\vskip 12pt
\centerline{Figure 1.}

The following two corollaries are
sample applications of Theorem 1. They control convergence of
hyperbolic alternating link complements in the geometric
topology. We will show that the only limit points are
the `obvious' ones, namely augmented alternating link
complements, as defined by Adams in [2].

\noindent {\bf Corollary 2.} {\sl A complete finite volume hyperbolic 3-manifold is
the limit of a sequence of distinct hyperbolic alternating link 
complements if and only if it is a hyperbolic augmented alternating
link complement.}

\noindent {\bf Corollary 3.} {\sl The set of all hyperbolic alternating
and augmented alternating link complements is a closed
subset of the set of all complete finite volume hyperbolic 3-manifolds, in the
geometric topology.}

The upper bound on volume is proved by using techniques
related to the Gromov norm [4].
We will show that the volume of $S^3 - K$ is at most
the volume of a link complement
with a diagram having $4t(D)$ crossings. By constructing
an explicit ideal triangulation for this link complement,
we find an upper bound for its volume.

The lower bound is established by using a theorem of
Agol [3]. When a finite volume hyperbolic 3-manifold $M$
contains a properly embedded 2-sided incompressible
boundary-incompressible surface $S$, Agol established a lower bound on the
volume of $M$ in terms of the `guts' of $M - {\rm int}({\cal N}(S))$.
In our case, $M$ is the complement of $K$, and $S$
is the orientable double cover of one of the two
`checkerboard' surfaces arising from an alternating diagram.

\vskip 18pt
\centerline {\caps 2. The upper bound on volume}
\vskip 6pt

We will use the fact [8] that if a compact orientable
hyperbolic 3-manifold $M$ is obtained by Dehn filling
another hyperbolic 3-manifold $N$, then the volume
of $M$ is less than the volume of $N$. The 3-manifold $N$ 
we will use is the exterior
of the link $J$ that is obtained by replacing each twist
of the diagram $D$ with a tangle containing at most
six crossings. This tangle is composed of the
two original strings of the twist, but with all but two
(respectively, all but one) of its crossings removed,
depending on whether the twist contained an even
(respectively, odd) number of crossings. Those
two strings are then encircled with a simple closed
curve, as in Figure 2, known as a {\sl crossing circle}. 
(There is one exception to this:
if two of these crossing circles cobound an annulus in the complement
of the remaining link components, then only one of these should
be used.) The resulting link $J$ is
an augmented alternating link, and hence is hyperbolic [2].
The link $K$ is obtained from $J$ by performing $1/q$ surgeries, 
for certain integer values of $q$, on the crossing circles. Hence,
${\rm Volume}(S^3 - K) < {\rm Volume}(S^3 - J)$.
If we alter the diagram of $J$ near each crossing circle
by removing the residual crossing(s) of the twist, 
the result is a new link $L$.
By [1], ${\rm Volume}(S^3 - L) = {\rm Volume}(S^3 - J)$.

\vskip 18pt
\centerline{\psfig{figure=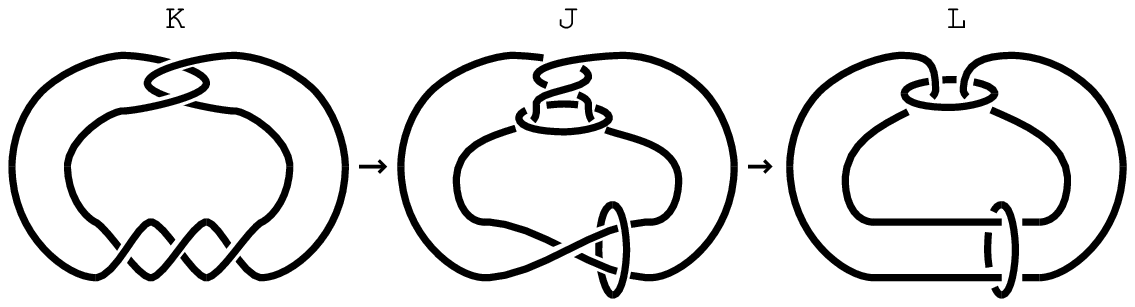}}
\vskip 6pt
\centerline{Figure 2.}

Note the diagram $D_L$ of $L$ is prime and connected, 
since $D$ is prime, connected and not the standard diagram of the $(2,n)$-torus link. 
Hence, it determines a decomposition of $S^3 - L$
into two ideal polyhedra with their faces identified in
pairs [4]. Here, we are using the term {\sl polyhedron} in
quite a general sense: a 3-ball with a connected graph
in its boundary that contains no loops and no valence one
vertices. An {\sl ideal polyhedron} is a polyhedron with its
vertices removed. The edges of this ideal polyhedral decomposition
of $S^3 - L$ are
vertical arcs, one at each crossing. The faces are the regions
of the diagram, twisted a little near the crossings so that
their boundaries run along the link and the edges, and
so that the interiors of the faces are disjoint. The remainder
of $S^3$ is two open 3-balls, which we take to be the interiors
of the two ideal polyhedra $P_1$ and $P_2$.

Note that the intersection of the 2-skeleton with the
boundary tori of the link exterior is a 4-valent
graph. Using Euler characteristic, the number of
complementary regions of this graph is equal to the number
of vertices. The former is the number of vertices of
the ideal polyhedra. The latter is $2 c(D_L)$, where
$c(D_L)$ is the number of crossings in the diagram $D_L$, since
there are two vertices for each edge of the polyhedral
decomposition.

We now subdivide the faces of the polyhedra 
with more than three boundary edges
into triangles by coning from an ideal vertex. We wish
to calculate the resulting number of triangles.
Let $V$, $E$ and $F$ be the total number of ideal vertices,
edges and faces in the boundary of the two ideal
polyhedra. So, 
$$V - E + F = \chi(\partial P_1) + \chi(\partial P_2) = 4.$$ 
The number of triangles is the sum, over all faces
of the polyhedra,
of the number of sides of the face minus two. This is
$2E - 2F = 2V - 8 = 4 c(D_L) - 8$. 

Now collapse
each bigon face of the polyhedra to a single edge. Some
care is required here, since it is {\sl a priori} possible
that there is a cycle of bigons, glued together along their
edges. However, an examination of the diagram $D_L$
gives that the bigons are in fact disjoint. 

In each polyhedron, there is a vertex with valence at
least four. Otherwise, the boundary graph is a single triangle
or the tetrahedral graph, and it is straightforward to check
that these graphs do not arise. For example, observe that each region
of $D_L$ has an even number of sides and there must be more than two
non-bigon regions. So each polyhedron ends up with more than four
triangular faces. Triangulate each polyhedron by coning
from this vertex. The result is an ideal
triangulation of the complement of $L$
with at most $4c(D_L) - 16$ tetrahedra. 

This allows us to bound the volume of the complement
of $L$. We homotope each ideal tetrahedron so that it lifts
to a straight simplex in the universal cover ${\Bbb H}^3$.
First homotope each edge, preserving its ends, so that
it is either a geodesic or has been entirely homotoped into
a cusp. Then homotope each ideal 2-simplex so that it
is straight, but possibly degenerate. Then do the
same for the 3-simplices. The volume of each resulting ideal
3-simplex is at most $v_3$. Hence, ${\rm Volume}(S^3 - L)
\leq v_3(4c(D_L) - 16) = v_3 (16 t(D) - 16)$.
This proves the right-hand inequality of the main theorem.
Note again that we did not use that $D$ is alternating.

\vskip 18pt
\centerline {\caps 3. The lower bound on volume}
\vskip 6pt

The lower bound on the volume of a hyperbolic alternating
link is proved using the following theorem of Agol [3]. 
It deals with a finite volume hyperbolic 3-manifold $M$
containing a properly embedded incompressible 
boundary-incompressible surface $S$.
We denote $M - {\rm int}({\cal N}(S))$ by $M_S$.

\noindent {\bf Theorem.} (Agol [3]) {\sl Let $M$ be an
orientable hyperbolic 3-manifold containing a properly
embedded orientable boundary-incompressible incompressible surface $S$. Then 
$${\rm Volume}(M) \geq -2 \, v_3 \, \chi({\rm Guts}(M_S)).$$

}

We refer the reader to [3] for a full description of 
the `guts' terminology. Essentially, the pair
$(M_S, \partial {\cal N}(S) \cap \partial M_S)$ 
has an associated characteristic submanifold $\Sigma$, which
is a canonical collection of $I$-bundles and Seifert fibred spaces
embedded in $M_S$, and the {\sl guts} of $M_S$ is the closure of 
the complement of $\Sigma$. We refer to $P =
\partial M \cap \partial M_S$ as the {\sl parabolic locus}.
It is a collection of annuli and tori.

Note that the assumption that $S$
is orientable can be dropped, providing the surface
$\tilde S ={\rm cl}(\partial {\cal N}(S) - \partial M)$ is
incompressible and boundary-incompressible. For if we apply Agol's theorem to 
$\tilde S$, then $M_{\tilde S}$ is a copy of ${\cal N}(S)$
and a copy of $M_S$. The former is an $I$-bundle and hence
a component of the characteristic submanifold of $M_{\tilde S}$.
Hence, ${\rm Guts}(M_{\tilde S}) = {\rm Guts}(M_S)$.
 
In our case, $M$ is the exterior of the alternating link $K$,
and $S$ is one of the two {\sl checkerboard} surfaces $B$ and $W$,
arising from a diagram of $K$. These surfaces arise by colouring the regions
of the diagram black and white, so that regions meeting along
an arc of the link projection have different colours.
If all the faces with the same colour are glued together,
twisted near the crossings, the result is one of the checkerboard surfaces.
However, instead of using
the given diagram $D$ of $K$, it is convenient to work with a diagram
that is in addition {\sl twist-reduced}.
This means that whenever a simple closed
curve in the diagram intersects the link projection transversely
in four points disjoint from the crossings, and two of these points are adjacent to
some crossing, and the remaining two points are adjacent to
some other crossing, then this curve bounds a subdiagram
that consists of a (possibly empty) collection of bigons 
arranged in a row between these two crossings. An equivalent
pictorial definition of a twist-reduced diagram is
given in Figure 3.

\vskip 18pt
\centerline{\psfig{figure=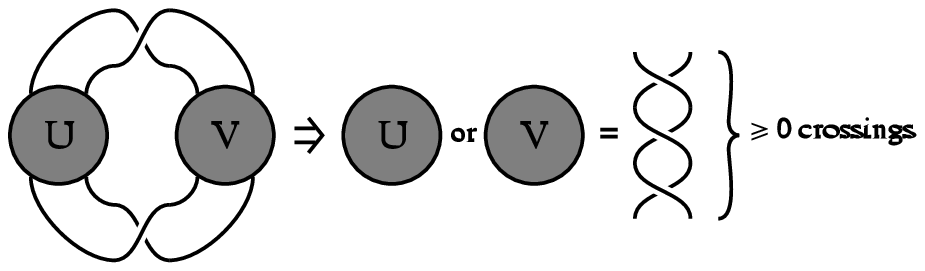}}
\vskip 6pt
\centerline{Figure 3.}

Any prime alternating link has a twist-reduced prime
alternating diagram. For if a diagram is not twist-reduced,
it decomposes as in Figure 3, where the top and bottom
crossings lie in different twists. There is a sequence of
flypes that amalgamates these into a single twist,
reducing the twist number of the diagram. However,
we need to know more than this. We will prove the following in
\S4.

\noindent {\bf Lemma 4.} {\sl Let $K$ be a link with a connected prime
alternating diagram $D$. Then $K$ has a connected prime alternating
twist-reduced diagram $D'$ with at least $(t(D)/2+1)$ twists.}

We let $B$ and $W$ be the black and white checkerboard surfaces
for the twist-reduced diagram $D'$. Let $r_B(D')$ and $r_W(D')$
be the number of black and white non-bigon regions of $D'$.
We will prove the following theorem in \S5.

\vfill\eject
\noindent {\bf Theorem 5.} {\sl Let $D'$ be a prime alternating
twist-reduced diagram of $K$, let $M$ be the exterior of $K$,
and let $B$ and $W$ be the checkerboard surfaces for $D'$.
Then
$$\eqalign{
\chi({\rm Guts}(M_B)) &= 2 - r_W(D')\cr
\chi({\rm Guts}(M_W)) &= 2 - r_B(D'). \cr}$$
}

Note that the diagram $D'$ induces a planar graph,
with a vertex at each twist and an edge for each
edge of $D'$ that is not adjacent to a bigon region.
Denote the number of its vertices, edges and faces
by $V$, $E$ and $F$. Then $2E = 4V$, since it is 4-valent.
Hence,
$$2 = V - E + F = -V + F = -t(D') + r_B(D') + r_W(D').$$
The lower bound on volume follows rapidly from these results.
Adding the inequalities in Agol's theorem applied
to $B$ and $W$, we obtain
$$\eqalign{
{\rm Volume}(S^3 - K) &\geq
-v_3
\big(
\chi({\rm Guts}(M_B)) + \chi({\rm Guts}(M_W)) \big) \cr
&= -v_3(4 - r_B(D') - r_W(D')) \cr
&= v_3(t(D') - 2) \cr
&\geq v_3(t(D)/2 - 1).}$$

\vskip 12pt
\centerline{\psfig{figure=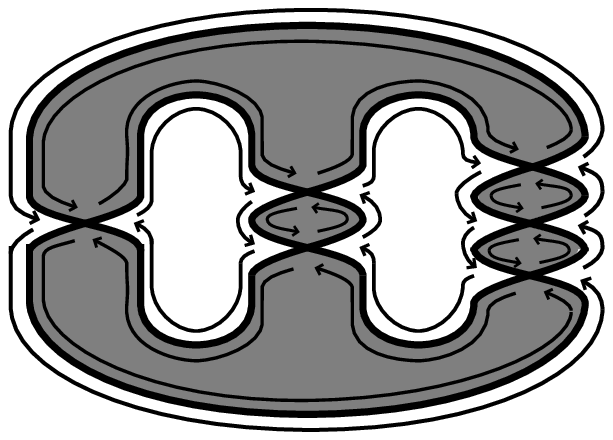,width=2in}}
\vskip 12pt
\centerline{Figure 4.}

The proof of Theorem 5 relies heavily on
the ideal polyhedral decomposition of a link complement
arising from a connected prime diagram, as described in \S2.
Its 2-skeleton is the union of the two checkerboard
surfaces. When the diagram is
alternating, the boundary graphs of the two polyhedra $P_1$
and $P_2$ are particularly simple. They are just copies of the
underlying 4-valent graph of the link projection [4].
Each region of one polyhedron is glued to the
corresponding region of the other, with a rotation
that notches the face around one slot in either
a clockwise or anti-clockwise direction, depending
on whether the region is coloured white or black.
Thurston compared this gluing procedure to the gears of
a machine [8].

Throughout much of this paper we will consider a surface
$S$ properly embedded in $M_B$ (or $M_W$). It will
intersect the parabolic locus $P = \partial M \cap
\partial M_B$ in a (possibly empty) collection of
transverse arcs. It will
be incompressible and also {\sl parabolically incompressible}
(see Figure 5) which means that there is no embedded disc $E$ in $M_B$
such that
\item{$\bullet$} $E \cap S$ is a single arc in $\partial E$;
\item{$\bullet$} the remainder of $\partial E$ is an arc in $\partial M_B$
which has endpoints disjoint from $P$ and which
intersects $P$ in at most one transverse arc;
\item{$\bullet$} $E \cap S$ is not parallel in $S$ to an arc in $\partial S$ 
that contains at most one component of $\partial S \cap P$.

\noindent Also no component of $S$ will be a boundary-parallel disc 
or a 2-sphere.

\vskip 18pt
\centerline{\psfig{figure=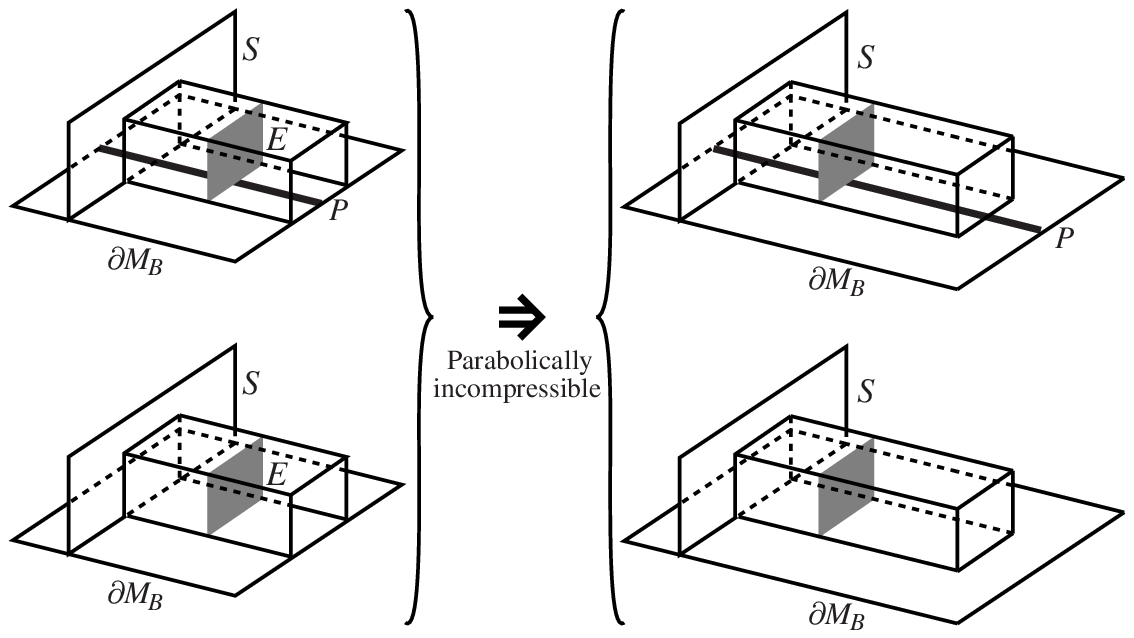}}
\vskip 6pt
\centerline{Figure 5.}

A fairly standard argument gives that such a surface $S$ can be ambient isotoped
(leaving $P$ invariant) into {\sl normal form}, which means it satisfies the following
conditions:
\item{$\bullet$} the intersection of $S$ with each of the
ideal polyhedra $P_i$ is a collection of properly embedded discs;
\item{$\bullet$} each such disc intersects any boundary edge
of $P_i$ at most once;
\item{$\bullet$} the boundary of each such disc cannot enter
and leave an ideal vertex through the same face of $P_i$;
\item{$\bullet$} $S$ intersects any face of $P_i$ in a collection of
arcs, rather than simple closed curves;
\item{$\bullet$} no such arc can have endpoints in the same
ideal vertex of $P_i$, or in a vertex and an edge that are adjacent;
\item{$\bullet$} no component of intersection between $S$
and $\partial P_i$ forms the boundary of a regular neighbourhood
of an edge.

In fact, for this argument to work, we need to know that
$\partial {\cal N}(B) \cap M_B$ is incompressible and
boundary-incompressible in $M_B$, but we will check this below. 
The last of the above conditions is non-standard. It can
be guaranteed since $S$ is properly embedded in $M_B$ rather
than $M$. 

Using normal surfaces, we can prove the following lemma,
which we need in order to apply Agol's theorem. This is
a stronger version of a result of Menasco and Thistlethwaite
(Proposition 2.3 of [6]) that asserts that $B$ and $W$
are incompressible in $M$.

\noindent {\bf Lemma 6.} {\sl The surfaces $\tilde B = \partial {\cal N}(B)
\cap M_B$ and $\tilde W = \partial {\cal N}(W) \cap M_W$
are incompressible and boundary-incompressible in $M$.}

\noindent {\sl Proof.} If there is a compression
disc for $\tilde B$, say, then there is 
one $S$ in normal form. The intersection of $S$ with
$W$ can contain no simple closed curves, since
such a curve would lie in a face of the polyhedral decomposition. 
So it is a collection of arcs. Suppose initially that there
is at at least one such arc. An outermost one in $S$
separates off a disc $S'$ in $S$ that lies in some $P_i$.
Its boundary intersects the edges of $P_i$
twice and misses the vertices. The boundary graph
of $P_i$ is a copy of the link diagram. Since the
diagram is prime, we deduce that $\partial S'$ intersects the
same edge of $P_i$ twice. This contradicts the definition
of normality. Therefore, $S \cap W$ is empty. But, then
$S$ lies entirely in one $P_i$, with its boundary
in a black face. It is therefore not a compression disc. 

If there is a boundary-compression disc for $\tilde B$,
then there is one in normal form. As above, we may assume
that $S \cap W$ is empty. But $S$ then lies in one $P_i$, with
$\partial S$ running over a single ideal vertex and avoiding
all edges of $P_i$. This contradicts the fact that the
diagram is prime. $\square$

\vskip 18pt
\centerline {\caps 4. The characteristic decomposition of a link diagram}
\vskip 6pt

Analogous to the characteristic submanifold of a 3-manifold,
in this section we define the characteristic decomposition of
a connected prime link diagram $D$.
We will consider simple closed curves in the diagram.
Always these will be disjoint from the crossings and
will intersect the link projection transversely.
We will, at various points, isotope these curves.
Always, the isotopy will leave the crossings fixed
throughout. In this section, we will make this proviso
without further mention.

A {\sl square} is a simple closed curve that intersects
the link projection four times. The link projection
divides it into four arcs, which we call its
{\sl sides}. A square is {\sl essential} if
it is not homotopically trivial in the complement
of the crossings, or, equivalently, it intersects four distinct 
edges of the link projection. A square is 
{\sl characteristic} if it is essential, does not
separate off a single crossing, and any other square
can be isotoped off it. Taking one isotopy class
of each characteristic square, and isotoping them
so that they are all disjoint, gives a collection of
squares which we term the {\sl characteristic collection}.

\vskip 18pt
\centerline{\psfig{figure=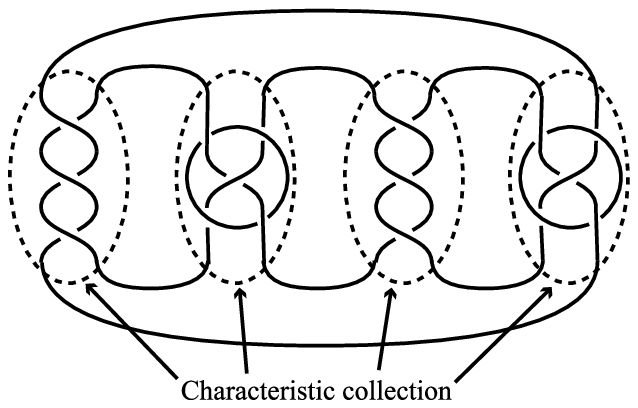}}
\vskip 12pt
\centerline{Figure 6.}
\vfill\eject

Consider a collection ${\cal C}$ of disjoint essential squares. Define
a connected complementary region $R$ of ${\cal C}$ to be a {\sl product region}
if each component of
intersection between $R$ and the black regions (or each component
of intersection between $R$ and the white regions) intersects
${\cal C}$ and the crossings in a total of two components.
(See Figure 7.) 
Thus, $R$ is a copy of a 2-sphere with some open discs removed,
the open discs and the crossings lying on the equator of the 2-sphere,
and the black (or white) regions forming a regular neighbourhood
of the remainder of the equator. The reason for the name
product region is as follows. If we view this 2-sphere as the
boundary of a 3-ball $B$, then $B$ minus an open regular
neighbourhood of the crossings is
the product of a closed interval and a disc. 
The white regions in $R$ form the
horizontal boundary of this product, and the black regions
lie in the vertical boundary.

\vskip 18pt
\centerline{\psfig{figure=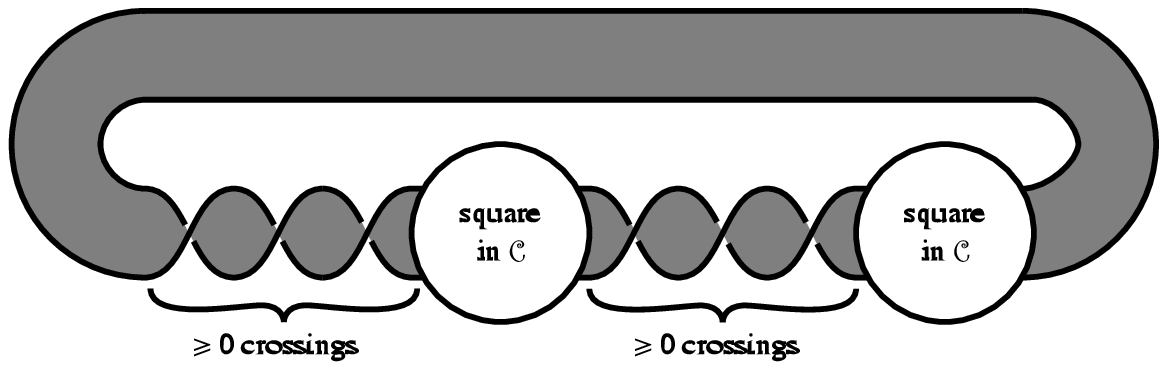}}
\vskip 12pt
\centerline{Figure 7.}

\noindent {\bf Lemma 7.} {\sl Let $S$ and $S'$ be
essential squares, isotoped in the complement of
the crossings to minimise $|S \cap S'|$. Then
they intersect in zero or two points,
and in the latter case, the two points of intersection
lie in distinct regions of the diagram with the
same colour.}

\noindent {\sl Proof.} The four sides of $S$ run through
distinct regions of the diagram, as do the four sides
of $S'$. A side of $S$ intersects a side of $S'$
at most once. So $|S \cap S'|$ is at most four.
However, if it is exactly four, then $S$ and $S'$
run through the same regions, hence are isotopic
and so can be ambient isotoped off each other.
So, $|S \cap S'|$ is at most three. It is even and therefore
either zero or two. In the latter case, the regions
containing these two points cannot have opposite
colour, as each arc of $S - S'$ and $S' - S$ would
then intersect the link projection an odd number of
times. So one of the four complementary regions of
$S \cup S'$ would have only two points of intersection
with the link projection in its boundary. As the diagram
is prime, it would contain a single arc of the link
projection and no crossings, and hence $S$ and $S'$
could be ambient isotoped off each other. $\square$

\noindent {\bf Lemma 8.} {\sl Let $R$ be a product
complementary region of a non-empty collection ${\cal C}$ of
disjoint essential squares. Then any essential square in
$R$ that is characteristic must be parallel to square in
${\cal C}$.}

\noindent {\sl Proof.} Suppose that the white regions in $R$ form
the horizontal boundary of $R$. Pick one such region $E$.
Then (isotopy classes of) essential squares in $R$ are in
one-one correspondence with (isotopy classes of) properly
embedded arcs in $E$ with endpoints in distinct black regions.
Let $S$ be an essential square in $R$ that is not parallel
to any curve in ${\cal C}$ and that does not enclose
a single crossing. Let $\alpha$ be the corresponding
arc in $E$. Then each component of $\partial E - \partial \alpha$
contains an arc of intersection with a black region
that is disjoint from $\partial \alpha$. Pick an arc
$\beta$ in $E$ joining these black regions. This 
corresponds to an essential square that cannot be 
isotoped off $S$. Hence, $S$ is not characteristic. $\square$

\noindent {\bf Lemma 9.} {\sl Let ${\cal C}$ be a collection
of disjoint non-parallel essential squares, such that
\item{1.} any essential square in the complement of
${\cal C}$ lies in a product complementary region, 
or is parallel to a square of ${\cal C}$, or encloses
a single crossing;
\item{2.} if two product complementary regions of ${\cal C}$
are adjacent, they have incompatible product structures;
\item{3.} no square of ${\cal C}$ encloses a single crossing.

\noindent Then ${\cal C}$ is characteristic.}

\noindent {\sl Proof.} We will first show that each
square in ${\cal C}$ is characteristic. So, consider
an essential square $S$ and isotope it to intersect 
${\cal C}$ minimally. Suppose that $S$ is not disjoint
from ${\cal C}$, intersecting some curve $S'$ of
${\cal C}$. By Lemma 7, $S'$ must intersect $S$
twice, and these two points of intersection lie
in different regions of the diagram with the
same colour, white say. Since all the curves
in ${\cal C}$ are disjoint, any other curve of ${\cal C}$
intersecting $S$ must do so in the white regions.
Hence, two components of $S - {\cal C}$ are arcs
intersecting the link projection twice. Consider
one such arc $\alpha$, and let $S_1$ be the curve
of ${\cal C}$ containing its endpoints. Let $\beta$
be the arc(s) of $S - {\cal C}$ adjacent to $\alpha$,
and let $S_2$ be the component(s) of ${\cal C}$
touching $\beta$, one of which is $S_1$.
Then ${\cal N}(S_1 \cup \alpha)$ is a product
region. By (1), each boundary curve either is
parallel to a curve of ${\cal C}$, or lies in a product
complementary region of ${\cal C}$, or encircles a 
single crossing. Hence,
${\cal N}(S_1 \cup \alpha)$ lies in a product
complementary region $R_1$ of ${\cal C}$
(apart from an annular strip running along
$S_1$). Similarly, ${\cal N}(S_2 \cup \beta)$
lies in a product complementary region $R_2$ of ${\cal C}$.
They have compatible product structures, contradicting
(2). Hence, each square in ${\cal C}$ is
characteristic.

Now, any characteristic square can be isotoped off
${\cal C}$, and so, by (1), is either parallel
to a curve in ${\cal C}$ or lies in a product
complementary region. In the latter case, it
must be parallel to a curve in ${\cal C}$, by 
Lemma 8. So, ${\cal C}$ consists of all the characteristic
squares. $\square$

\noindent {\bf Lemma 10.} {\sl The characteristic
collection satisfies (1), (2) and (3) of Lemma 9.}

\noindent {\sl Proof.} Pick a maximal collection
of disjoint non-parallel essential squares.
This satisfies (1).
Remove all squares that encircle a single crossing.
The resulting collection satisfies (1) and (3).
Remove a square if the two complementary regions
either side of it are product regions and
their union is again a product region. The
collection still satisfies (1) and (3). Repeat this
process as far as possible. The resulting collection
satisfies (1), (2) and (3), and hence is the characteristic
collection by Lemma 9. $\square$

\noindent {\bf Lemma 11.} {\sl Let ${\cal C}$ be
a collection of disjoint non-parallel essential squares that
bound a product region $R$. Suppose also that ${\cal C}$
is not a single square enclosing a single crossing.
Then $R$ extends to a product complementary region of 
the characteristic collection.}

\noindent {\sl Proof.} Extend ${\cal C}$ to maximal
collection of disjoint non-parallel essential squares that 
are disjoint from the interior of $R$. Then,
as in the proof of Lemma 10, reduce it to collection
satisfying (1), (2) and (3). The resulting
collection is characteristic by Lemma 9 and, by construction,
$R$ extends to a product complementary region of
the characteristic collection. $\square$

\noindent {\bf Corollary 12.} {\sl Each twist with more than
one crossing lies in a product complementary region
of the characteristic collection.}

\noindent {\sl Proof.} Apply Lemma 11 to a curve encircling
the twist. $\square$

We are now in a position to prove Lemma 4.

\noindent {\sl Proof of Lemma 4.} Let ${\cal C}$ be the characteristic
collection of squares for the diagram $D$. Each complementary product region
is of the form shown in Figure 7.
If this contains at least two twists, then there is a sequence of
flypes that amalgamates all these twists into one. Perform all
these flypes for all product regions, giving a new diagram $D'$.
The squares ${\cal C}$ in $D$ give squares ${\cal C}'$ in
$D'$. Note that a complementary region of ${\cal C}$ is a product
region before flyping if and only if the corresponding
region of ${\cal C}'$ is a product after flyping.
Adjacent product regions of ${\cal C}$ have incompatible
product structures, by Lemma 10, and so adjacent product
regions of ${\cal C}'$ have incompatible product structures.
So, by Lemma 9, ${\cal C}'$ is characteristic.
Note that, by construction, each product complementary
region of ${\cal C}'$ contains at most one twist.

We claim that $D'$ is twist-reduced. If not, it has a decomposition
as in Figure 3, where neither of the tangles $U$ nor $V$ is a row
of bigons as shown. The squares $\partial U$ and $\partial V$
are essential and cobound a product region. By Lemma 11,
this is part of a product region of ${\cal C}'$. This 
contains at most one twist, and therefore either $U$
or $V$ is a row of bigons, as shown in Figure 3.
Therefore $D'$ is twist-reduced.

We now have to prove that $t(D') \geq t(D)/2+1$. 
Let ${\cal C}_+$ be the characteristic collection
${\cal C}$ of $D$, together with a curve enclosing
each crossing that is not part of a longer twist.
By Corollary 12, each twist lies in a product complementary
region of ${\cal C}_+$.

Define a graph $G$ having a vertex for each complementary
region of ${\cal C}_+$.
Two vertices are joined by an edge if and only if
the corresponding regions are adjacent. Denote the
vertex set by $V(G)$ and the subset of the vertices
that arise from a product region with at least one
twist by $T(G)$. For any vertex $v$ of $V(G)$, let
$\lambda(v)$ be its valence. Note that the
vertices of valence one correspond to innermost
regions of the diagram, which are necessarily
a single twist, and hence lie in $T(G)$. Note also that the valence of a 
vertex in $T(G)$ is at least the number of twists
that the corresponding region of $D$ contains. Now,
$G$ is a tree, and so, by Euler characteristic,
$$\sum_{v \in V(G)}(\lambda(v) - 2) = -2.$$
Therefore,
$$\sum_{v \in T(G)}(\lambda(v) - 2) = -2 - \sum_{v \in V(G) - T(G)}
(\lambda(v) - 2) \leq -2.$$
So,
$$
t(D) \leq \sum_{v \in T(G)} \lambda(v)
= \sum_{v \in T(G)}(\lambda(v) - 2) + 2|T(G)| \leq - 2 + 2|T(G)| = - 2 +
2t(D').$$
Rearranging this gives that $t(D') \geq t(D)/2 + 1$. $\square$

\vskip 18pt
\centerline{\caps 5. The guts of the checkerboard surfaces' exteriors}
\vskip 6pt

Let $D'$ be a twist-reduced prime alternating diagram.
Let $B$ be its black checkerboard surface and let $R$ be
its white bigon regions. Then we may view $R$ as
a collection of discs properly embedded in $M_B$. Each disc 
of $R$ intersects the parabolic locus $P = \partial M \cap M_B$ 
transversely twice. Hence, ${\cal N}(R \cup P)$
is an $I$-bundle $E$ embedded in $M_B$. It is
part of the characteristic submanifold of $M_B$.
So, the guts of $M_B$ is the guts of ${\rm cl}(M_B - E)$, 
where the latter is given parabolic locus ${\rm cl}(\partial E - \partial M_B)$.

We can identify ${\rm cl}(M_B - E)$
explicitly. If each twist of $D'$ consisting of white
bigons is removed and replaced with a single crossing,
the result is an alternating diagram $\hat D$ (of
a new link). Let $\hat B$ be its black
checkerboard surface. Then 
${\rm cl}(M_B - E)$ is homeomorphic
to $S^3 - {\rm int}({\cal N}(\hat B))$, the
homeomorphism taking parabolic locus to parabolic
locus. Hence, it suffices to analyse
$\hat B$ and $\hat D$. It is clear that
$\hat D$ is prime, but it may not be twist-reduced.
Nevertheless, it is {\sl black-twist-reduced},
which means that the implication of Figure 3 holds
whenever the regions of the figure adjacent to both
$U$ and $V$ are black. For example, the diagram in
Figure 8 is black-twist-reduced, but not reduced.
It in fact arises from the diagram in Figure 6 by 
assigning a checkerboard colouring, and then removing
the white bigons.

\noindent {\bf Theorem 13.} {\sl Let $D$ be
a black-twist-reduced prime alternating
diagram, with no white bigon regions.
Let $B$ be its black checkerboard
surface. Then 
$$\chi({\rm Guts}(M_B)) = \chi(M_B).$$}

\vskip 18pt
\centerline{\psfig{figure=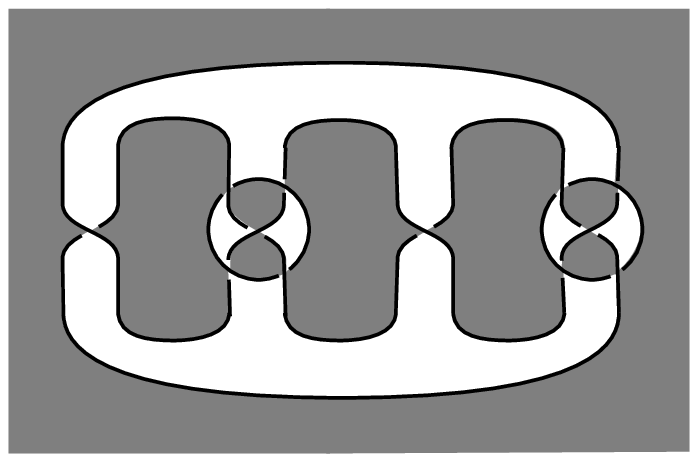}}
\vskip 12pt
\centerline{Figure 8.}

Theorem 5 is an immediate corollary, since we may
apply Theorem 13 to $\hat D$ and deduce that
$\chi({\rm Guts}(M_B)) = \chi(S^3 - {\rm int}({\cal N}(\hat B))$.
But $S^3 - {\rm int}({\cal N}(\hat B))$ is a regular neighbourhood of a graph,
with two vertices and $r_W(D')$ edges. One vertex lies above the
diagram, and one vertex lies below, and each edge runs
between the vertices through a non-bigon white region of $D'$.
So its Euler characteristic is $2 - r_W(D')$.

We now embark on the proof of Theorem 13. So, suppose that
$D$ is a black-twist-reduced prime alternating
diagram with no white bigon regions. Let $A$ be a 
characteristic annulus for the characteristic
submanifold of $M_B$. It is incompressible,  and
disjoint from the parabolic locus.
It may be parabolically compressible, in which
case it parabolically compresses to a {\sl product
disc}, which is a disc properly embedded in 
$M_B$ intersecting the parabolic locus in two
transverse arcs. This product disc is {\sl essential} in the sense
that it is not boundary-parallel, provided $A$ is
not boundary-parallel. We will prove the following result.

\noindent {\bf Theorem 14.} {\sl Let $D$ be a connected 
black-twist-reduced
prime alternating diagram with no white bigon regions.
Let $B$ be its black checkerboard surface. Then
\item{(i)} $M_B$ contains no essential product discs, and
\item{(ii)} any incompressible parabolically-incompressible
annulus $A$ that is disjoint from the parabolic locus 
separates off a Seifert fibred solid torus
subset of the characteristic submanifold of $M_B$.

}

\noindent {\sl Proof.}
Suppose first that there is an essential product disc $S$.
It is incompressible and parabolically incompressible,
and so can be ambient isotoped into normal form. The intersection
of $S$ with $W$ is a collection of arcs. Each of these
arcs must run between distinct components of
$S \cap \partial {\cal N}(B)$. For, otherwise,
we may consider an arc of $S \cap W$ with endpoints in the
same component of $S \cap \partial {\cal N}(B)$, that is outermost in $S$.
This separates off a subdisc of $S$ that lies in a single
ideal polyhedron and that has two points of intersection with
the polyhedron's edges. Since the link diagram is prime, these
must lie in the same edge, contradicting normality.

Suppose now that $S \cap W$ is non-empty.
An outermost arc of $S \cap W$ in $S$ separates
off a disc $S_1$. This lies in one of the ideal
polyhedra of the link complement, and so gives a
curve in the link diagram, as shown in Figure 9. Let
$N$ be the region of the diagram containing $S_1 \cap W$.
Note that $S_1 \cap N$ cannot separate off a single crossing
in $\partial N$. For $D$ would then decompose as in Figure 3.
However, the fact that $D$ is black-twist-reduced would imply
that one of $U$ or $V$ in Figure 3 would contain a row of
white bigons, whereas $D$ contains no white bigons.

\vskip 18pt
\centerline{\psfig{figure=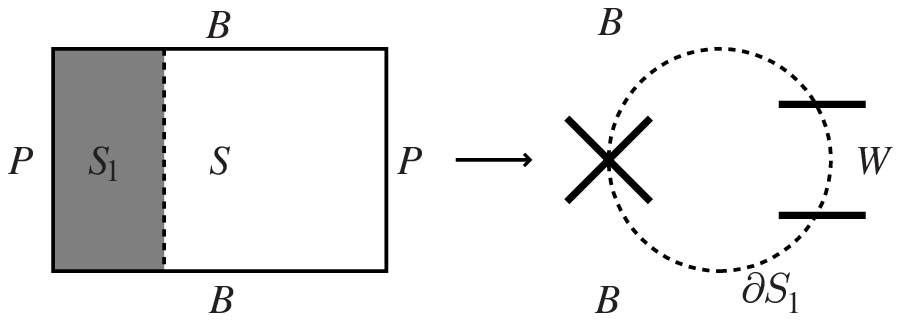}}
\vskip 6pt
\centerline{Figure 9.}

The arc $S_1 \cap W$ is part of a normal disc 
$S_2$ in the other ideal polyhedron. Both ideal polyhedra have 
the same boundary graph. So we may super-impose 
$\partial S_1$ and $\partial S_2$. Perform
a small ambient isotopy so that they miss the crossings,
ensuring that the new white side of $\partial S_1$
is disjoint from $\partial S_2$. The result is two
essential squares. Perform a further isotopy in the
complement of the crossings so that they have minimal
intersection. Since $S_1 \cap N$ cannot separate off a single
crossing in $\partial N$, $\partial S_1$ and $\partial S_2$
must intersect in $N$. So, in the white
regions, they have only a single point of intersection.
(See Figure 10.) This contradicts Lemma 7.

\vskip 18pt
\centerline{\psfig{figure=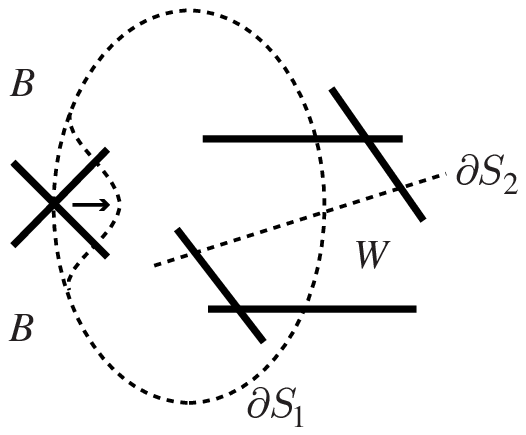}}
\vskip 6pt
\centerline{Figure 10.}

Hence, the assumption that $S \cap W$ is non-empty is impossible.
The disc $S$ therefore lies entirely above or below the diagram.
It runs over two crossings. The fact that $D$ is black-twist-reduced
gives that these two crossings are joined by a row
of white bigons, again a contradiction. Thus, $M_B$
contains no essential product discs.

Now consider an incompressible parabolically-incompressible
annulus $A$ that is disjoint from the parabolic locus.
Again, we may place $A$ in normal form. As in the
case of the product disc, the arcs $A \cap W$ must
run between distinct boundary components of $A$. 

Consider a disc $S_2$ lying between two adjacent arcs $\alpha_1$ and
$\alpha_3$ of $A \cap W$. Each is part of
discs $S_1$ and $S_3$ properly embedded in the other ideal polyhedron.
As before, super-impose $\partial S_2$ and $\partial S_i$
($i = 1$ or $3$), and minimise their intersection with an ambient isotopy.
If $R_i$ is the region of the diagram containing $\alpha_i$,
then $\alpha_i$ cannot separate off a single crossing
in $\partial R_i$. For this would imply the existence of
a parabolic compression disc. See Figure 11. 
Hence, $\partial S_2$ and $\partial S_i$ intersect in the
white region $R_i$. By Lemma 7, they have exactly one
other point of intersection, in some other white region.
Since $\partial S_1$ and $\partial S_3$ are disjoint,
the only possibility for $\partial S_1$, $\partial S_2$ 
and $\partial S_3$ is as shown in Figure 12. 

\vskip 18pt
\centerline{\psfig{figure=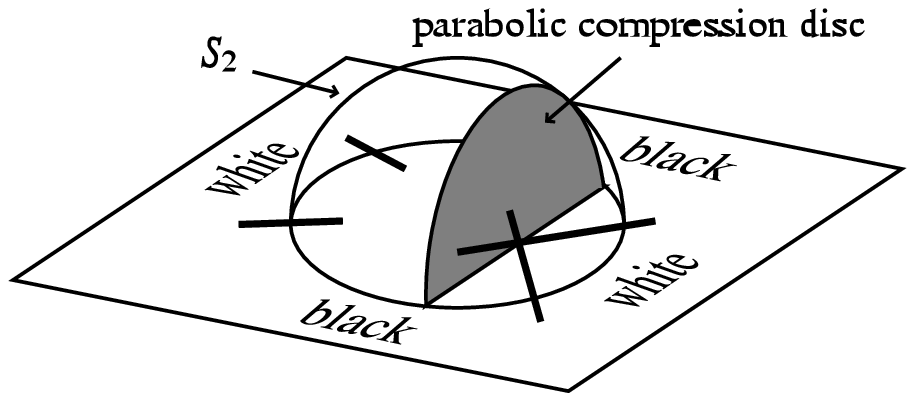}}
\vskip 6pt
\centerline{Figure 11.}

\vskip 24pt
\centerline{\psfig{figure=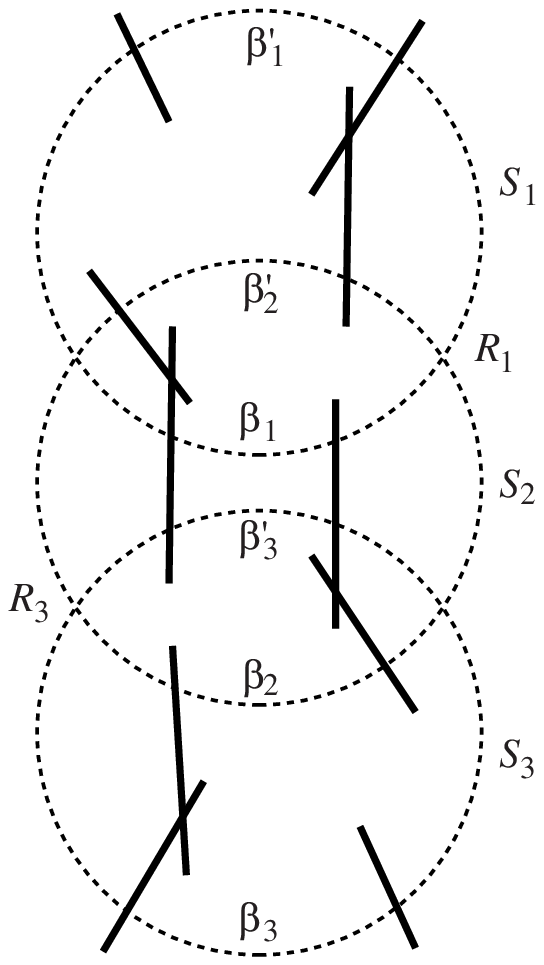,width=1.5in}}
\vskip 18pt
\centerline{Figure 12.}

We term a disc subset of the diagram a {\sl unit}
if its boundary is an essential square, and 
its intersection with some white region of the
diagram contains a single crossing. Note that
a unit is uniquely specified by its two
black sides, since these sides extend uniquely
to a square and a square cannot bound units on
both sides in a black-twist-reduced diagram with no
white bigons. Two units are {\sl fused} if
\item{$\bullet$} they are disjoint;
\item{$\bullet$} two of their black sides are parallel; and
\item{$\bullet$} either the white boundary sides that separate
off a crossing are not adjacent, or at least one of
the units contains a single crossing.

\noindent Note that, if a unit is fused to some other
unit, then either it contains a single crossing or
its boundary is characteristic. For if it does not contain
a single crossing, then its boundary has product
regions with incompatible product structures on 
both sides. By Lemma 11, these extend to product
complementary regions of the characteristic collection.
So, the boundary of the unit is characteristic.

\vskip 18pt
\centerline{\psfig{figure=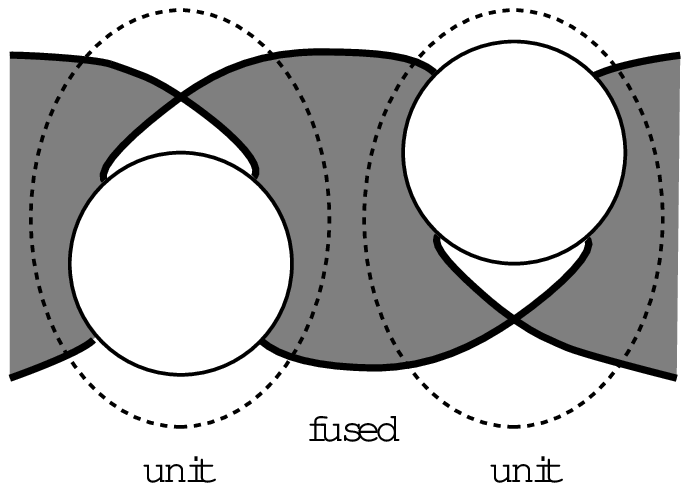}}
\vskip 6pt
\centerline{Figure 13.}

Denote the arcs $\partial A \cap S_i$ by $\beta_i$
and $\beta'_i$, where $\beta_1$, $\beta_2$ and $\beta_3$
all lie in the same component of $\partial A$.
We deduce from Figure 12 that $\beta_1$ and $\beta_2$ form
opposite black edges of a unit, as do $\beta_2$ and $\beta_3$,
and that these units are fused. This argument may be
applied to the two squares ($S_2$ and $S_4$) either side of
$S_3$. If $\beta_4$ is the arc of $\partial A \cap S_4$
adjacent to $\beta_3$, we deduce that $\beta_3$ and $\beta_4$
form opposite black edges of a unit, which is fused
to the unit bounded by $\beta_2$ and $\beta_3$. 
Continue in this way. The resulting units all
have non-intersecting boundaries (up to isotopy) since they are
characteristic or separate off a single crossing. They cannot be non-trivially nested,
since they are fused to other units. Hence, the
only possibility is that the units start to repeat.
That is, they are all fused in a circular fashion,
as in Figure 14.

We observe from Figure 12 that, for either of the polyhedra $P_i$,
each curve of $A \cap \partial P_i$ encircles
two units. These curves are all disjoint.
Also, if two discs of $A \cap P_i$
have a single disc between them (for example, $S_1$ and $S_3$), 
then two of their black edges are parallel. There is therefore
a constant number ($n$, say) of copies of each
curve. Thus, $A$ is $n$ parallel copies of a surface. Since $A$ is
connected, $n=1$. Therefore, the intersection
of $A$ with the boundary of each ideal polyhedron
is as shown in Figure 14. (There, the case
$|A \cap W| = 6$ is shown, but $|A \cap W|$ may
be any even integer greater than three.)

So, for each ideal polyhedron $P_i$,
one component of $P_i - A$ is a ball with boundary
a product region disjoint from the ideal vertices.
The intersection of the product region with $W$
is two discs, one in the central region of Figure 14,
one in the outer region. These two balls, one in each
ideal polyhedron, glue together to form the
required Seifert fibred solid torus. $\square$

\vskip 18pt
\centerline{\psfig{figure=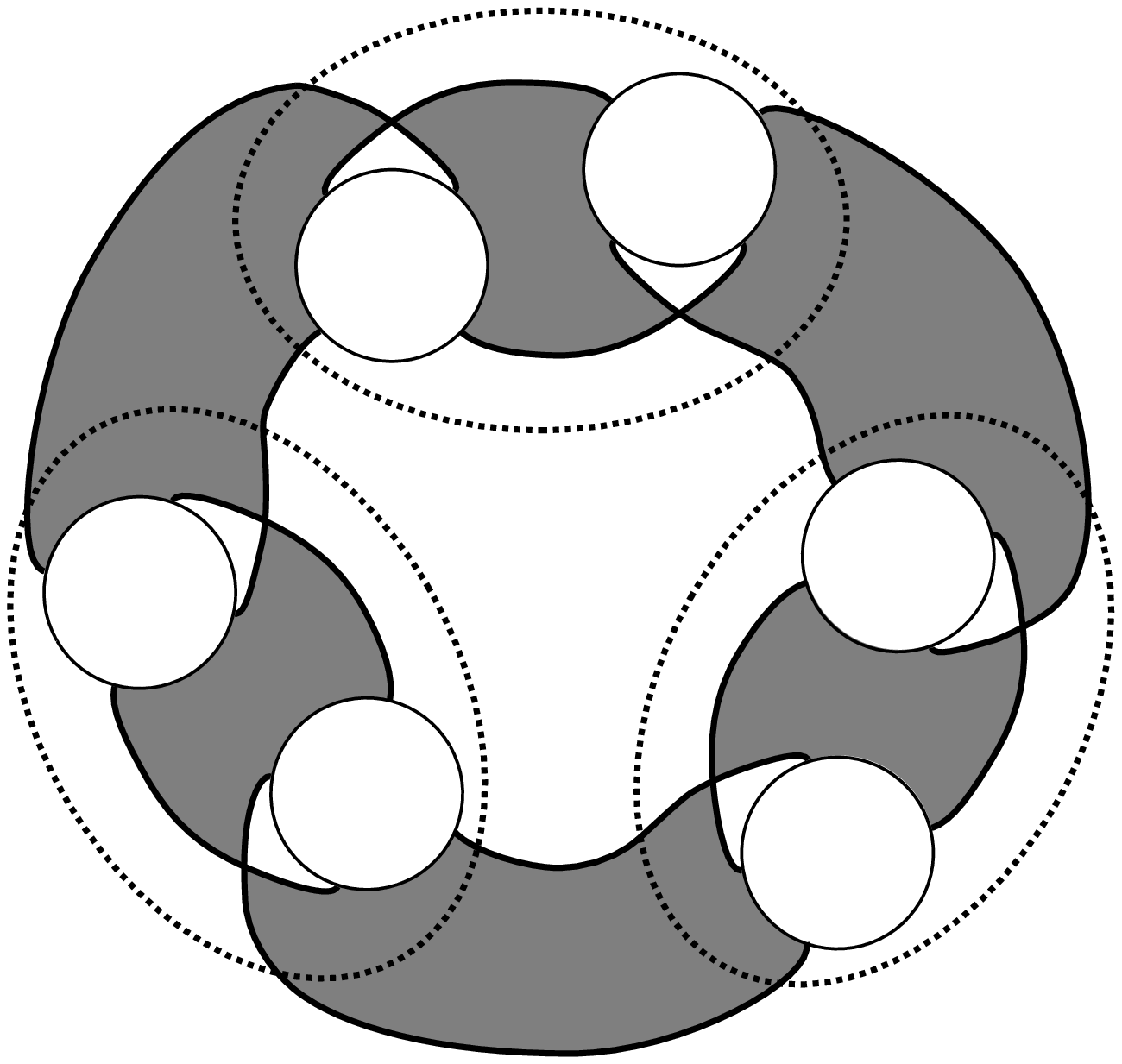,width=3in}}
\vskip 18pt
\centerline{Figure 14.}

\noindent {\sl Proof of Theorem 13.} Consider
a bounding annulus of the characteristic submanifold
of $M_B$. By Theorem 14, it is either boundary parallel
or separates off a Seifert fibred solid torus. This
Seifert fibred solid torus is part of the characteristic
submanifold. We claim that the part of $M_B$ on the other side of this
annulus cannot also be part of the characteristic submanifold.
If it were, it would have to be an $I$-bundle.
If this intersected the parabolic locus, it would contain an
essential product disc, contradicting Theorem 14.
Its bounding annuli must each separate off a Seifert
fibred solid torus by Theorem 14. But then $M_B$ contains
no parabolic locus, which is a contradiction, proving the
claim. So, the characteristic
submanifold of $M_B$ is a collection of Seifert fibred
solid tori attached to the guts of $M_B$ via annuli.
Since each solid torus and each annulus has zero Euler characteristic,
$\chi({\rm Guts}(M_B)) = \chi(M_B)$. $\square$

\vfill\eject
\centerline{\caps 6. Application: convergent sequences}
\centerline{\caps of alternating link complements}
\vskip 6pt

As a sample application of the main theorem of this
paper, we show that the only possible limits of sequences
of hyperbolic alternating link complements converging
in the geometric topology are the `obvious' ones.

\noindent {\bf Corollary 2.} {\sl A complete finite volume hyperbolic 3-manifold is
the limit of a sequence of distinct hyperbolic alternating link 
complements if and only if it is a hyperbolic augmented alternating
link complement.}

\noindent {\sl Proof.} It is well known that a sequence of complete
finite volume hyperbolic 3-manifolds $M_i$ converges in the geometric
topology to a finite volume hyperbolic 3-manifold $M_\infty$ 
if and only if all $M_i$ sufficiently far along
this sequence are obtained by Dehn filling $M_\infty$,
so that none of the surgery slopes have a constant subsequence. 
Hence, one way of obtaining a convergent sequence of alternating
hyperbolic link complements is to start with a single
augmented alternating link and perform surgeries along
the augmenting unknots, where the surgery coefficients
are chosen so that the resulting diagrams are alternating.
Thus, any hyperbolic augmented alternating link complement
is certainly the limit of a sequence of distinct hyperbolic alternating
link complements. We must show that the converse also holds.

Consider a collection of
distinct hyperbolic alternating knots $\{ K_i \}$ with complements
that converge in the geometric topology to a finite
volume hyperbolic 3-manifold. Their volumes
are bounded. Therefore, by the main theorem of this
paper, the twist numbers of their alternating diagrams
are bounded. In \S2, we
showed how a link with a diagram $D$ is obtained by
surgery along crossing circles in a link with a diagram
having at most $6t(D)$ crossings. When the original diagram
is alternating, the new diagram is augmented alternating.
Since the twist numbers of the $K_i$ are bounded, so
are the crossing numbers of these augmented alternating links $L_i$. 
Thus, there are only finitely many possible $L_i$. We may therefore
pass to a subsequence in which $L_i$ is a constant link. 
The $K_i$ in this subsequence have alternating diagrams
obtained by surgically replacing each crossing circle tangle with a twist.
Pass to a subsequence so that, at every twist, the number of 
crossings is either constant or tends to infinity.
Then, the links in this subsequence are all obtained
by Dehn filling a single augmented alternating link $L$,
so that none of the surgery coefficients have a constant
subsequence. So, $L$ is the limit of the sequence. $\square$

\noindent {\bf Corollary 3.} {\sl The set of all hyperbolic alternating
and augmented alternating link complements is a closed
subset of the set of all complete finite volume hyperbolic 3-manifolds, in the
geometric topology.}

\noindent {\sl Proof.} Consider a convergent
sequence of hyperbolic alternating and augmented alternating
link complements. We wish to show that the limit hyperbolic
manifold is an alternating or augmented alternating link
complement. Arbitrarily close to each augmented
alternating link in the geometric topology, there is an alternating
link. Hence, we may assume that the sequence consists only
of alternating links. If it has a constant subsequence,
the limit is an alternating link complement. If it has no constant
subsequence, the limit is an augmented alternating
link complement, by Corollary 2. $\square$

\vskip 24pt

\+ Mathematical Institute,\cr
\+ Oxford University, \cr
\+ 24-29 St Giles',\cr
\+ Oxford OX1 3LB, \cr
\+ England.\cr

\vfill\eject
\centerline{\caps Appendix}
\centerline{\caps Ian Agol and Dylan Thurston}

In this appendix, we improve on Lackenby's upper bound 
on the volume of links in terms of the twist number of
the projection diagram, proven in section 2 of the paper. 
We will use the same notation as in section 2. 

\noindent {\bf Theorem.} {\sl 
Given a projection diagram $D$ of a link $K$ with twist number $t(D)$, 
${\rm Vol}(S^3-K)\leq 10 v_3(t(D)-1)$. Moreover, there is a 
sequence of links $K_i$ such that ${\rm Vol}(S^3-K_i)/t(D_i) \rightarrow
10 v_3$.}

\noindent {\sl Proof.} We use Lackenby's approach, by taking the link $K$ and
creating an augmented alternating link $L$ which has components
lying flat in the projection plane which are bound together
by crossing circles.

We will describe two decompositions of $S^3-L$, in
order to get an upper bound on the volume. The first decomposition
is obtained by taking the planar surface lying in the projection
plane, and taking the 2-punctured disks bounding each crossing
loop. The 2-punctured disks are totally geodesic
in the complement of $L$ [2], and there is a reflection through
the projection plane exchanging the polyhedra, and preserving
the planar surface in the projection plane. So this surface
is totally geodesic, and the 2-punctured disks are perpendicular
to the planar surface. 
These surfaces determine 
a decomposition of $S^3 - L$
into two ideal polyhedra $P_1$ and $P_2$ with their faces identified in
pairs (see Figure 15).
These polyhedra $P_1$ and $P_2$ are identical, and have
the property that the faces may be checkerboard colored dark
and white so
that the dark faces are triangles which come in pairs sharing
a vertex each, like a bow-tie. $S^3-L$ is obtained by
folding the bowties in $P_1$ along each vertex to glue the
pairs of triangles
together, then doubling along the rest of the faces. 

\vskip 18pt
\centerline{\psfig{figure=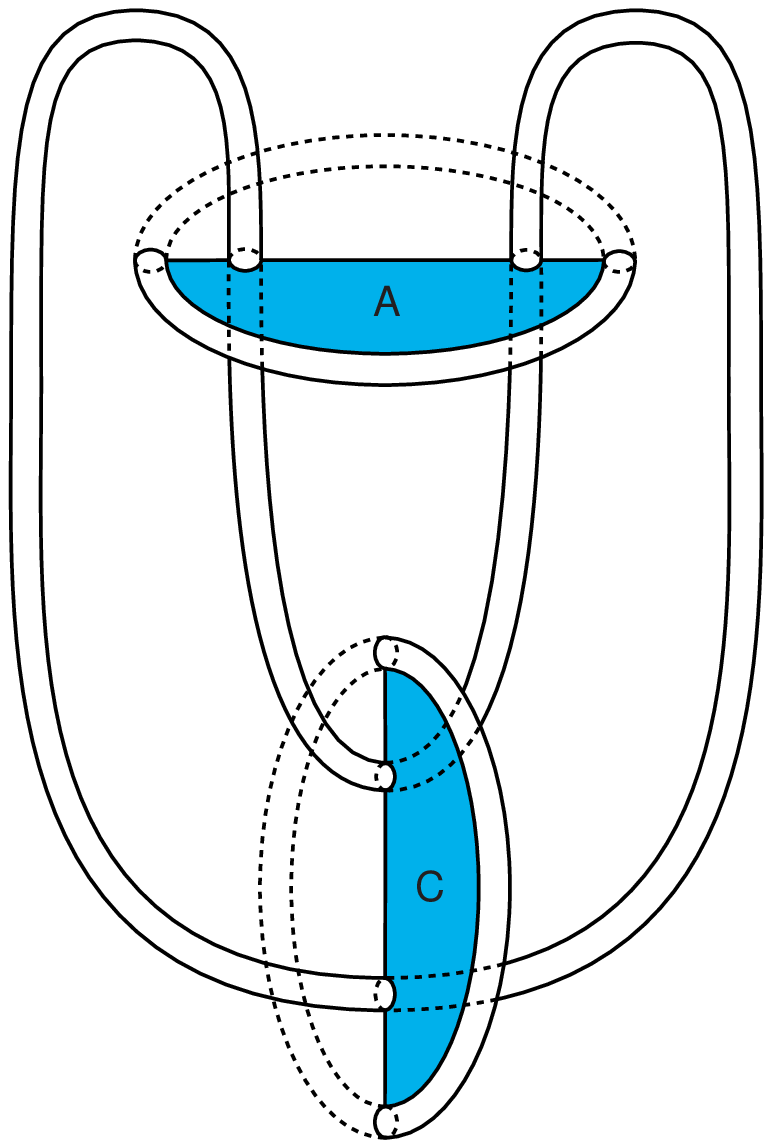,width=1.6in}}
\vskip 18pt
\centerline{\psfig{figure=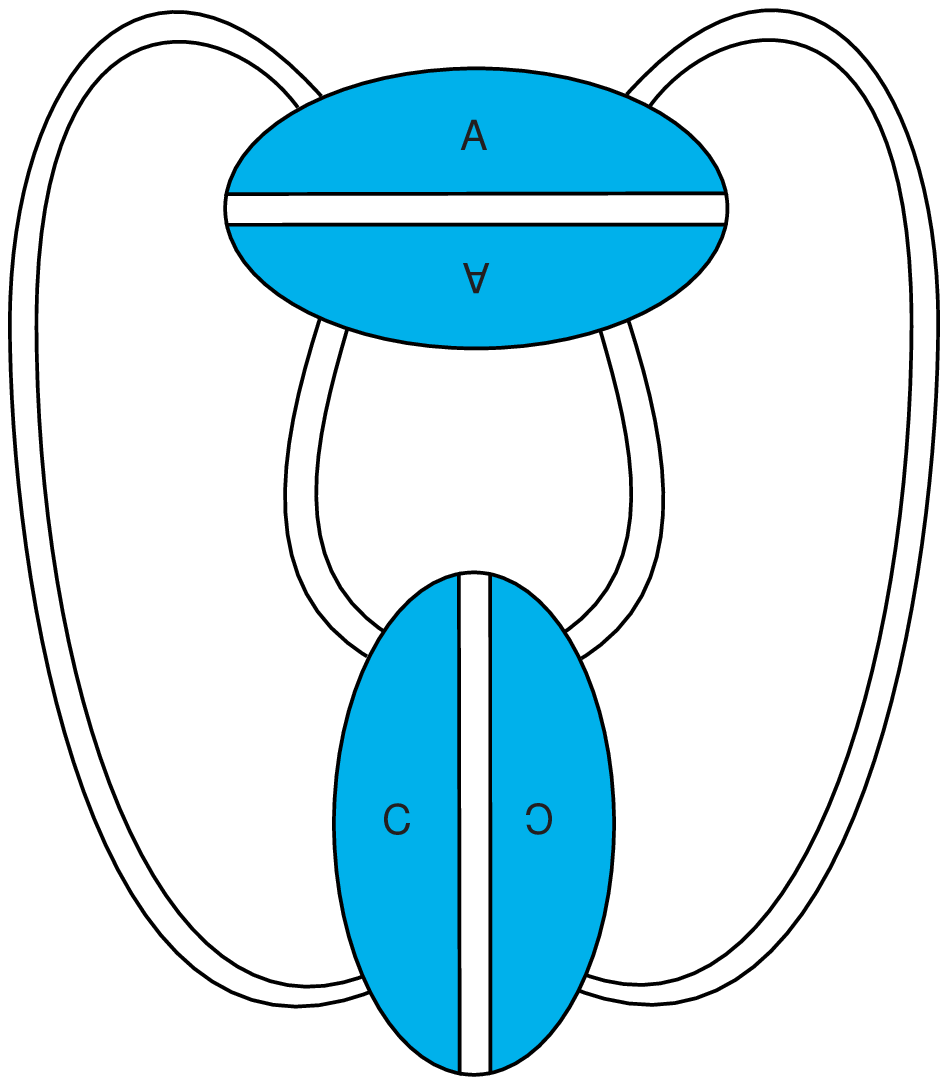,width=1.6in}}
\vskip 18pt
\centerline{\psfig{figure=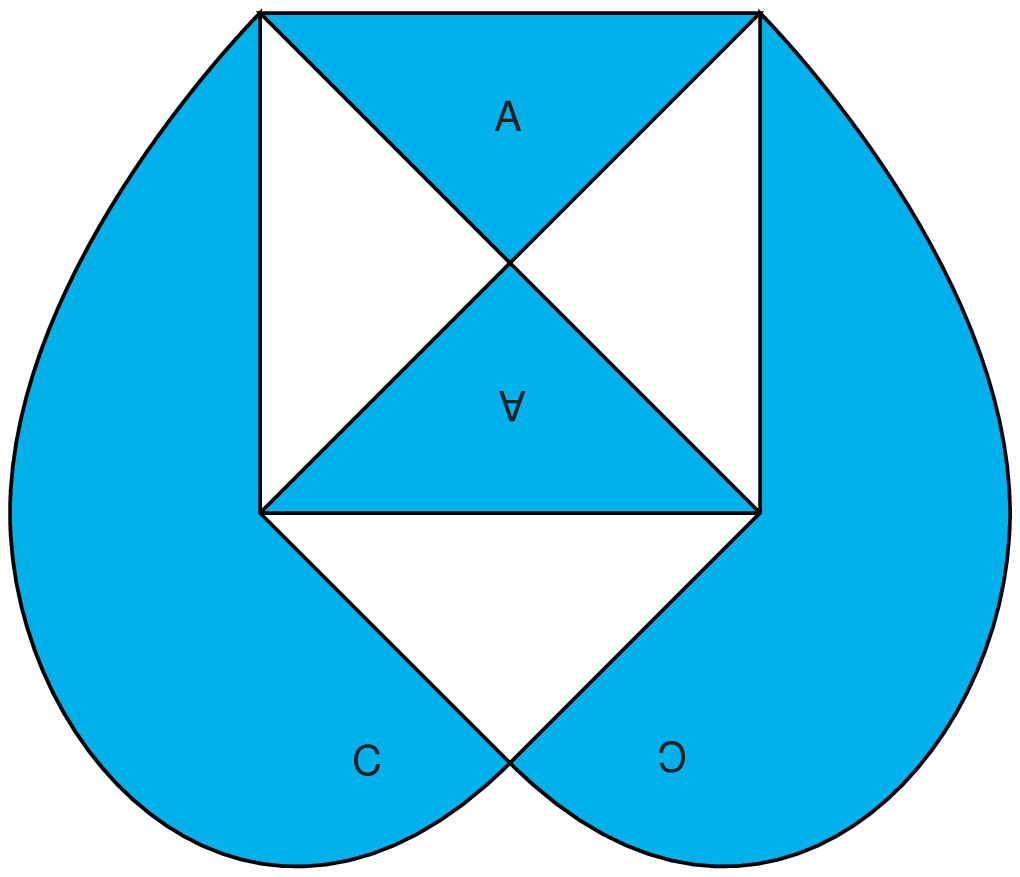,width=1.6in}}
\vskip 32pt
\centerline{Figure 15: Decomposing the complement of $L$ into ideal polyhedra}

The second decomposition is into tetrahedra. This is 
obtained by putting vertices
$v_1$ and $v_2$ in the
interior of $P_1$ and $P_2$, and coning the vertices to the faces
of the polyhedra. Each dark face of $P_1$ and $P_2$ gets coned
off to two tetrahedra, and each white face gets coned off to
two pyramids.
For each white face, we do a
stellar subdivision on the two pyramids containing it. That is,
we remove the face, and add an edge dual to the face connecting
$v_1$ and $v_2$. Then we add in triangles around the new edge
to divide the region into tetrahedra. If the face has $d$ edges,
then this divides the two cones into $d$ tetrahedra. 

To compute the total number of tetrahedra in this
triangulation, notice that each crossing loop contributes
6 edges to $P_1$ and $P_2$. Thus, the total number of 
edges in the white faces will be $6t(D)$, which 
contributes $6t(D)$ tetrahedra. Each crossing loop also
contributes 2 dark triangles each to $P_1$ and $P_2$, giving
4 tetrahedra when we cone off to $v_1$ and $v_2$.
Thus, we have a total of at most $10t(D)$ tetrahedra in this
decomposition. 

We may reduce the number of tetrahedra
by choosing an ideal vertex, and collapsing the edges
adjoining $v_1$ and
$v_2$ to this vertex. We then simplify the resulting
cell decomposition, by collapsing monogons to vertices,
bigons to single edges, and parallel triangles to single
triangles, to get an ideal triangulation. 
The vertex
we collapse to is adjacent to two dark faces and two white faces
(in the polyhedral decomposition into $P_1$ and $P_2$). Thus, when
we collapse the vertices $v_1$ and $v_2$ to an ideal
point, we collapse the 4 tetrahedra adjacent to the
dark faces to triangles.  All the white faces have degree $\geq 3$,
since we have assumed that no two crossing loops are parallel.
So we also collapse $\geq 6$ tetrahedra going through the two white
faces to triangles. 
Thus, we may find an ideal triangulation
with at most $10t(D)-10$ tetrahedra.

This allows us to bound the volume of  $S^3-L$, as
in section 2, by straightening the
triangulation, and bounding the volume of each
simplex by $v_3$, to conclude that ${\rm Vol}(S^3-L)\leq 10v_3(t(D)-1)$. 

Now, we consider the second claim of the theorem, which
shows that we have obtained the optimal constant. 
As motivation, we will consider the infinite augmented
link $C$ in ${\Bbb R}^3$  resembling a chain link fence,
which realizes exactly the upper bound on volume density
of $10v_3$ per crossing loop.

\vskip 18pt
\centerline{\psfig{figure=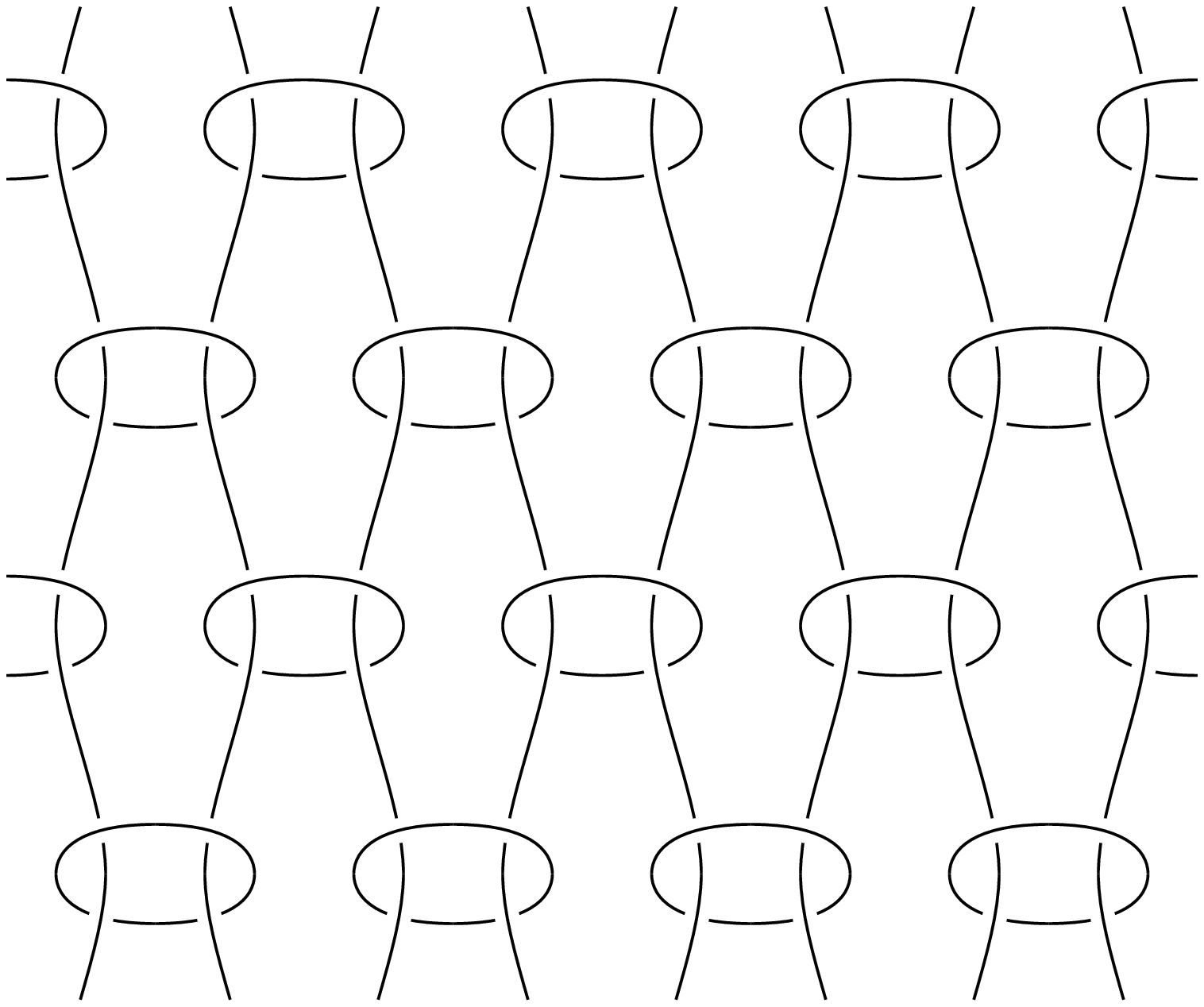,width=3in}}
\vskip 40pt
\centerline{Figure 16: The chain fence link}
 
Using the decomposition of the complement of $C$ 
into $P_1$ and $P_2$, we get
an infinite polyhedron with faces alternating between
triangles and hexagons, like a union of stars of David (see
Figure 17). 

This polyhedron has a natural realization as a right angled
polyhedron,
by taking the tessellation with regular
triangles and hexagons, and putting a circle around each
face. This is the same as taking the regular hexagonal
packing of circles, and putting a circle around each 
interstitial region (see Figure 17). 
So all the circles are either disjoint,
tangent, or orthogonal. Each circle bounds a geodesic
plane in the upper half space model of ${\Bbb H}^3$,
together cutting out a right-angled polyhedron.
Thus, we get a complete hyperbolic structure on the
link complement when we glue the polyhedra together to
get ${\Bbb R}^3-C$. 

\vskip 18pt
\centerline{\psfig{figure=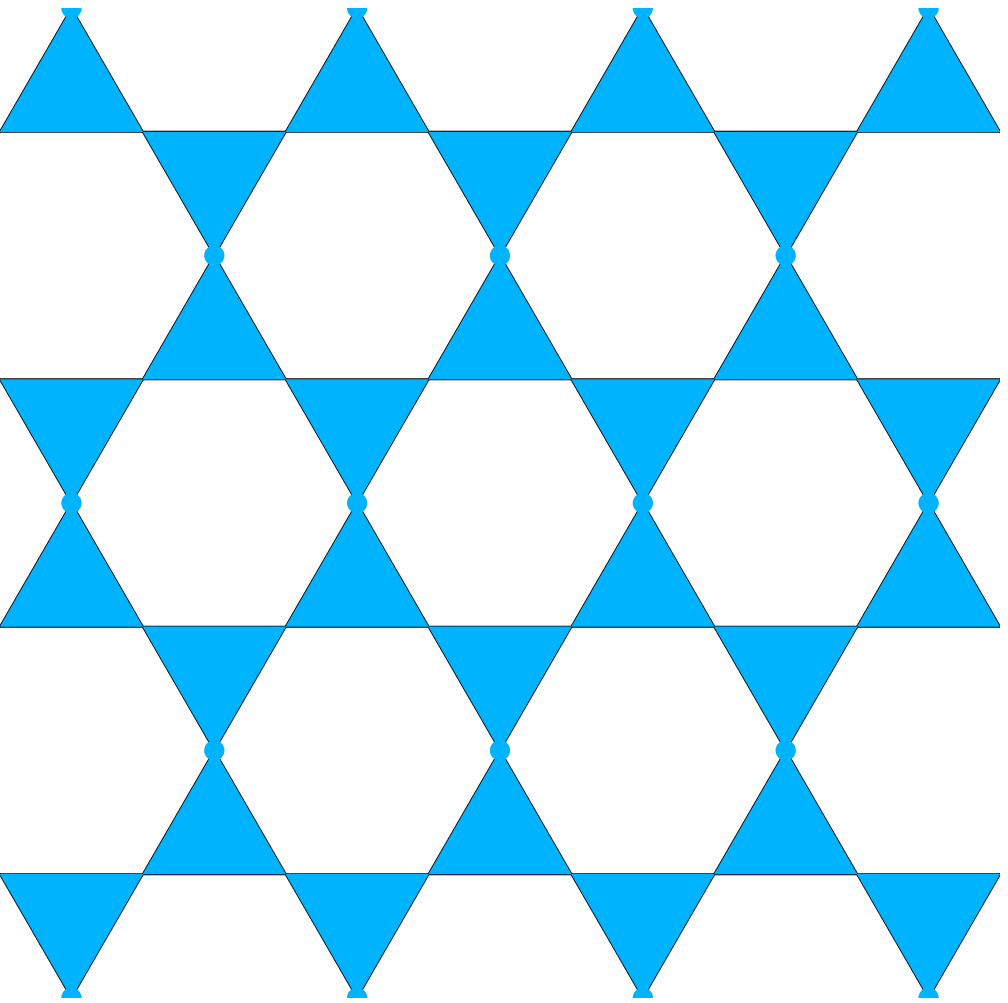,width=2.5in}}
\vskip 24pt
\centerline{\psfig{figure=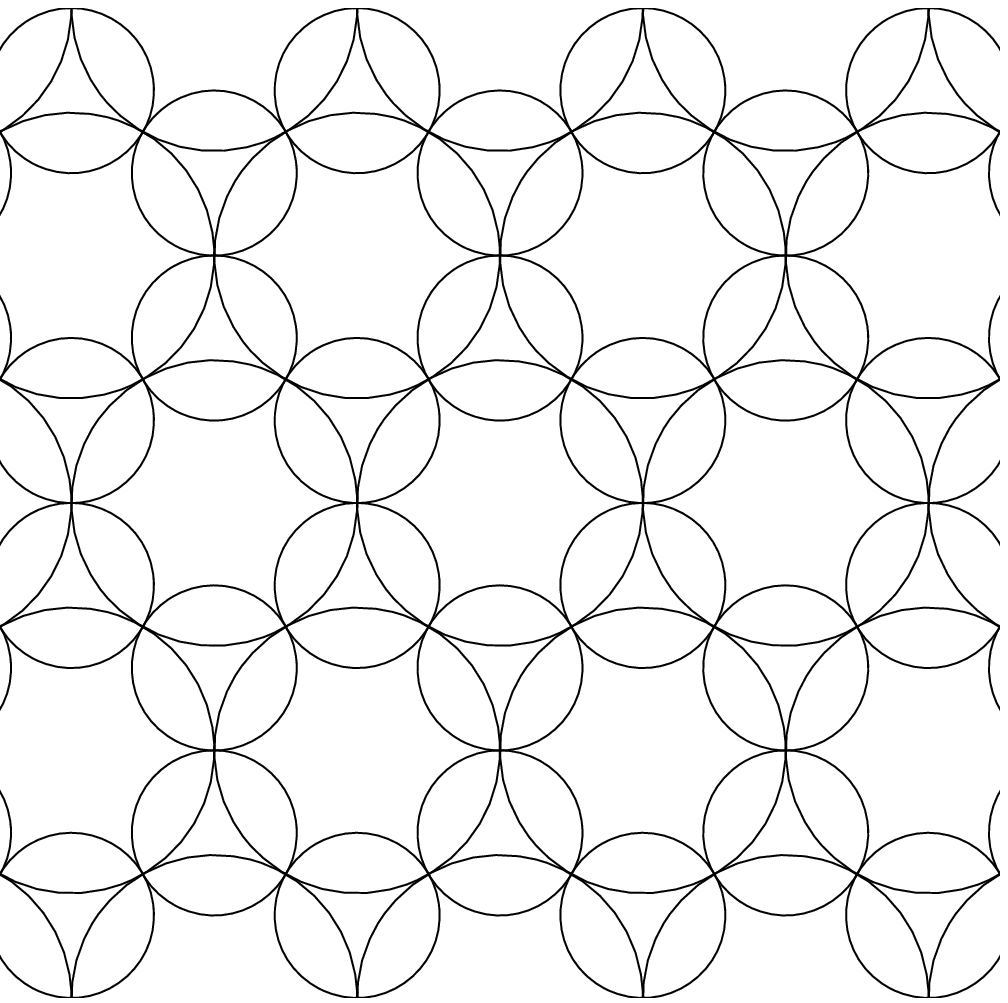,width=2.5in}}
\vskip 40pt
\centerline{Figure 17: The infinite polyhedron}
 
\noindent {\bf Remark.} There are actually many infinite augmented
alternating link complements with the same volume density. 
Take any decomposition of the triangles in the star of 
David tessellation (Figure 17)
into bow-ties. Then folding the bow-ties to glue the
triangles together and doubling
gives an ideal cell structure of an infinite augmented 
alternating link complement with volume density
$10v_3$ per crossing circle.

We may form a sequence of augmented  links $C_i$ by
taking the 2222-orbifold quotient of the link $C$, for
an increasing sequence of fundamental domains for the
2222-orbifold, and deleting the orbifold cone axes. 
The polyhedra of the links $C_i$ are obtained by taking
the 2222-orbifold quotient of the star of David tessellation,
and deleting the orbifold cone points. That these polyhedra are 
realized by right-angled polyhedra in ${\Bbb H}^3$ 
is a consequence of Andreev's theorem, see for example [8]. 
Coning the points
to a vertex, we get a division of this polyhedron into 
tetrahedra and cones on hexagons. Taking the circle
packing corresponding to the white faces of the polyhedron,
we get a circle packing on $S^2$ such that all but four
circles are adjacent to exactly 6 other circles.
Moreover, as $i \rightarrow \infty$, the majority of the circles
will have a packing around them which is combinatorially
equivalent to $n$ generations of the regular hexagonal
circle packing, where we may assume that $n \rightarrow \infty$
as $i \rightarrow \infty$. Sending the cone vertex to $\infty$ in
the upper half space model of ${\Bbb H}^3$, we see a
circle surrounded by $n$ generations of the regular 
hexagonal packing. A theorem of Rodin and Sullivan [7] then
shows that as $n\rightarrow \infty$, the ratios of radii of the
circles adjacent to the central  circle go to $1$. 
Thus, the majority of faces coned to $\infty$ then approaches either 
a regular tetrahedron, or a cone on a regular ideal 
hexagon. So as $i\rightarrow\infty$, the density of volume per
crossing circle of $S^3-C_i$ approaches
that of ${\Bbb R}^3-C$, which is $10v_3$. 
Then we may form links $K_i$ by taking $\pm 1/q$ surgeries
on the links $C_i$ to get alternating links, where 
we choose $q\rightarrow \infty$ fast enough that ${\rm Vol}(K_i)/
{\rm Vol}(C_i)\to 1$,
which we may do by Thurston's hyperbolic Dehn surgery
theorem [8]. $\square$

\vskip 18pt
\centerline {\caps References}
\vskip 6pt

\item{1.} {\caps C. Adams}, {\sl Thrice-punctured spheres in hyperbolic 
3-manifolds}, Trans. Am. Math. Soc. {\bf 287} (1985) 645-656.
\item{2.} {\caps C. Adams}, {\sl Augmented alternating link complements are
hyperbolic}, Low-dimensional Topology and Kleinian groups, London Math.
Soc. Lecture Note Series. {\bf 112}, Cambridge Univ. Press (1986)
\item{3.} {\caps I. Agol}, {\sl Lower bounds on volumes of hyperbolic Haken
3-manifolds}, Preprint (1999)
\item{4.} {\caps R. Benedetti and C. Petronio}, 
{\sl Lectures on Hyperbolic Geometry}, \hfill\break
Springer-Verlag (1992).
\item{5.} {\caps W. Menasco}, {\sl Closed incompressible surfaces in
alternating knot and link complements}, Topology {\bf 23}
(1984) 37-44.
\item{6.} {\caps W. Menasco and M. Thistlethwaite},
{\sl The classification of alternating links.}, Ann. Math.  {\bf 138}
(1993) 113--171.
\item{7.} {\caps B. Rodin and D. Sullivan}, 
{\sl The convergence of circle packings to the Riemann mapping.}
J. Differential Geom. {\bf 26} (1987) 349--360.
\item{8.} {\caps W. Thurston}, {\sl The geometry and topology of 
three-manifolds}, Lecture notes from Princeton University (1978--80).

\vskip 24pt

\+Department of Mathematics and Statistics, \cr
\+University of Melbourne, \cr
\+Parkville, VIC 30120, \cr
\+Australia.\cr

\vskip 12pt
\+Department of Mathematics,\cr
\+Harvard University,\cr
\+1 Oxford Street,\cr
\+Cambridge, MA 02138,\cr
\+USA.\cr

\end